# Notes on a paper of Tyagi and Holm:'A new integral representation for the Riemann Zeta function'

by

Michael Milgram[*], Consulting Physicist, Geometrics Unlimited, Ltd.
Box 1484, Deep River, Ont., Canada. K0J 1P0.

**Abstract:** It is shown that a new series representation of Riemann's Zeta function obtained by Tyagi and Holm leads to an interesting new recursion for Bernoulli numbers of even index as well as new representations of, and infinite series involving, Zeta functions of special (integer) argument.

**Results**

In a recent paper, Tyagi and Holm[1] present a new integral representation for the Riemann Zeta function. It is interesting that this representation yields some further new results for classical quantities. In particular, a new recursion formula for the Bernoulli numbers of even index emerges naturally.

From their new integral representation, Tyagi and Holm obtain a new series representation for the Riemann Zeta function $\zeta(s)$

$$\zeta(s) = \frac{\pi^{s-1}}{1-2^{1-s}} \sin(\frac{\pi s}{2}) \sum_{n=1}^{\infty} \frac{(2-2^{s-2n})}{\Gamma(2n+2)} \Gamma(2n-s+1)\zeta(2n-s+1), \quad (1)$$

reminiscent of similar results that can be found in Srivastava and Choi[2], Section 3.1. From this series representation a number of special cases emerge.

Consider the case $s=1$, where the series is convergent; it is well-known that $\zeta(s)$ has a singularity at this point, with residue equal to unity. In (1), the denominator term vanishes as $1-2^{1-s} \approx \log(2)(s-1)$ near $s=1$, leaving

$$\sum_{n=1}^{\infty} \frac{(1-2^{(-2n)})\zeta(2n)}{n(1+2n)} = \log(2) \quad (2)$$

This is a variation of Eqs 536 and 558 of ref. 2, which never appear in this simple combination.

Consider the case $s = 2m$, $m \in N$ (positive integers). The sum **(1)** must be broken up into three parts, according to $1 \leq n \leq m-1$, $n = m$ or $n > m$. The third part (infinite series) vanishes, and after evaluating the appropriate limits, one finds a new representation for $\zeta(2m)$ in terms of a finite series involving $\zeta(-2n+1)$. Explicitly

---

[*] mike@geometrics-unlimited.com

$$\zeta(2m) = \frac{\pi^{2m}(-1)^m}{2(1-2^{1-2m})}\left(-\frac{1}{\Gamma(2m+2)} + \sum_{n=1}^{m-1}\frac{(2-2^{2m-2n})\zeta(2n-2m+1)}{\Gamma(2n+2)\Gamma(2m-2n)}\right) \quad (3)$$

Using well-known substitutions for the Zeta functions of positive and negative (integer) arguments (ref. 2, Eq. 2.3(10) and 2.3(18) ) in terms of Bernoulli numbers $B_{2n}$ and solving, **(3)** gives a new recursion formula for the Bernoulli number $B_{2m}$

$$B_{2m} = \frac{\Gamma(2m+1)2^{(-2m)}}{(1-2^{(1-2m)})}\sum_{n=1}^{m}\frac{B_{2m-2n}(2-2^{(2m-2n)})}{\Gamma(2m-2n+1)\Gamma(2n+2)}. \quad (4)$$

**(4)** can also be rewritten as

$$\sum_{n=0}^{m}\binom{2m}{2m-2n}\frac{B_{2n}}{(2n-2m-1)}\left(\frac{2-2^{2n}}{2-2^{2m}}\right) = 0 \quad (5)$$

which differs from Euler's classical result (ref. 2, Eq. 1.6(6) )

$$\sum_{n=0}^{2m}\binom{2m+1}{n}B_n = 0. \quad (6)$$

Note that Euler's recursion includes the term $B_1 = -1/2$ whereas **(4)** does not. For example, the new recursion with *m=6* gives

$$B_{12} = (-11242 B_{10} - 12573 B_8 - 4092 B_6 - 385 B_4 - 6B_2)/2047 + \frac{1}{53222}B_0 \quad (7)$$
$$= -691/2730$$

whereas Euler's recursion yields

$$B_{12} = -22 B_{10} - 99 B_8 - 132 B_6 - 55 B_4 - 6B_2 - \frac{1}{13}B_0 - B_1 = -691/2730. \quad (8)$$

Finally, consider the case *s=2m+1* in **(1)**. The infinite series must again be broken into three parts, and the appropriate limits evaluated. This time the cancellation of singularities takes place between the zero belonging to the Zeta function of negative even argument (or the numerical term $(2-2^{(2m-2n)}) = 0$ when n=m), and the gamma function $\Gamma(2n-2m)$. The eventual result is

$$\zeta(2m+1) = \frac{(-1)^m \pi^{2m}}{(1-2^{-2m})} \left( -\frac{\log(2)}{\Gamma(2m+2)} + \sum_{n=1}^{m-1} \frac{(2-2^{1+2m-2n})\zeta'(2n-2m)}{\Gamma(2n+2)\Gamma(2m-2n+1)} \right.$$
$$\left. + \sum_{n=1}^{\infty} \frac{(2-2^{1-2n})\Gamma(2n)\zeta(2n)}{\Gamma(2m+2n+2)} \right) \quad (9)$$

giving the Zeta function of odd argument in terms of its first derivative of even negative argument and a convergent infinite series. For example,

$$\zeta(3) = -\frac{4}{3}\pi^2 \left( -\frac{\log(2)}{6} + \sum_{n=1}^{\infty} \frac{(2-2^{1-2n})\Gamma(2n)\zeta(2n)}{\Gamma(4+2n)} \right) \quad (10)$$

and

$$\zeta(5) = \frac{2}{15}\pi^2 \zeta(3) + \frac{16}{15}\pi^4 \left( -\frac{\log(2)}{120} + \sum_{n=1}^{\infty} \frac{(2-2^{1-2n})\Gamma(2n)\zeta(2n)}{\Gamma(6+2n)} \right). \quad (11)$$

In general, using (ref. 2, Eq. 2.3(22)) (with integer $n$)

$$\zeta'(-2n) = (-1)^n \frac{\Gamma(2n+1)}{2(2\pi)^{2n}} \zeta(2n+1) \quad n > 0 \quad (12)$$

**(9)** becomes

$$\zeta(2m+1) = \frac{(-1)^m \pi^{2m}}{(1-2^{-2m})} \left( -\frac{\log(2)}{\Gamma(2m+2)} + \sum_{n=1}^{\infty} \frac{(2-2^{1-2n})\Gamma(2n)\zeta(2n)}{\Gamma(2m+2n+2)} \right)$$
$$+ \frac{1}{(1-2^{-2m})} \sum_{n=1}^{m-1} \frac{(2^{2n-2m}-1)(-\pi^2)^n \zeta(2m-2n+1)}{\Gamma(2n+2)} \quad (13)$$

Such results are related to series that appear in ref. 2. (e.g. Eqs. 3.4(593) and 3.1(8) ), although the specific forms obtained here appear to be new. Notice that **(13)** reduces to a tautology (1=1) in the limit as $m \to \infty$. By a simple reordering, **(13)** also gives the sum of an infinite series involving zeta functions of integer even order in terms of a finite series of zeta functions of integer odd order

$$\sum_{n=0}^{\infty} \frac{(2^{-2n-2}-1)\Gamma(2n+2)\zeta(2n+2)}{\Gamma(2m+2n+4)} \right) = -\frac{1}{2}\sum_{n=0}^{m-1} \frac{(2^{-2n-2}-1)(-1)^n \zeta(2n+3)}{\Gamma(2m-2n)\pi^{2(n+1)}} \quad (13)$$
$$-\frac{\log(2)}{2\Gamma(2m+2)}$$

---

[1] Tyagi S. and Holm, C., "A new integral representation for the Riemann Zeta Function", preprint available from http://arxiv.org/abs/math-ph/0703073 (2007).

[2] Srivastava, H.M. and Choi, J, "Series Associated with the Zeta and Related Functions", Springer (2001).